\documentclass[12pt]{article}

\usepackage{geometry}
\usepackage{graphicx,epstopdf}
\graphicspath{{figs/}}
\usepackage{amsmath,amsfonts,amssymb,amstext,amsthm}
\newcommand{\diag}{\mathop{diag}}
\newcommand{\supp}{\mathop{supp}}
\usepackage[colorlinks=true,citecolor=blue,linkcolor=blue]{hyperref}
\newcommand{\bu}{\mathbf{u}}
\newcommand{\bw}{\mathbf{w}}
\newcommand{\R}{{\mathbb R}}
\newcommand{\C}{{\mathbb C}}

\newcommand{\cM}{{\mathcal M}}
\newcommand{\zihat}{\widehat{z_i}}
\newcommand{\intd}{\:\mathrm{d}}
\newcommand{\be}{\begin{equation}}
\newcommand{\ee}{\end{equation}}

\title{Mean field approximation\\ in conformation dynamics}
\date{May 25, 2009}

\author{Gero Friesecke \and Oliver Junge \and P\'eter Koltai
\thanks{Center for Mathematics, Technische Universit\"at M\"unchen, D-85747 Garching, Germany}}

\begin{document}
\maketitle 

\begin{abstract}
We propose a new approach to the transfer operator based analysis of the conformation dynamics of molecules.  It is based on a statistical independence ansatz for the eigenfunctions of the operator related to a partitioning into subsystems. Numerical tests performed on small systems show excellent qualitative agreement between mean field and exact model, at greatly reduced computational cost.
\end{abstract}



\thispagestyle{plain}

\section{Introduction} Conformation transitions of molecules reflect the global spatial/temporal
behaviour of the system, and in particular occur at much slower timescales compared to the
elementary frequencies of the system. In a small peptide with a dozen atoms,
the typical scale difference is already a factor $\sim 10^2$; for folding transitions in
proteins, it ranges between $10^8$ and $10^{16}$, placing these events well beyond the timescales
accessible via direct trajectory simulation. 

The transfer operator approach, introduced into MD by Deuflhard et al.\ in their fundamental paper \cite{DeDeJu99a} (cf.\ also \cite{Sc99a,ScHuDe01a,De03a,DeSc04a,DeWe05a,HoDiLa07a,WeKuWa07a}), 
allows to access long time effects through short time simulations, at the expense of needing to
simulate `ensembles' of initial conditions, or mathematically: to compute the evolution of densities
on space. In the latter approach, dominant conformations and their transition rates
can be identified via the leading eigenvalues and eigenfunctions of a suitable transfer operator.
The catch is that to compute the latter, which are functions on phase resp.\ configuration
space, the number of computational degrees of freedom 
of a direct discretization grows exponentially in the number of atoms.  

To overcome this problem, our goal in this paper is
to propose a mean field method for computing eigenstates of transfer operators, whose relationship to the exact 
eigenvalue problem for the transfer operator is reminiscent of that of Hartree-Fock theory to the many-particle 
Schroedinger equation in quantum chemistry, which overcomes an analogous `curse of dimension' problem. 
The mean field model only has the dimensionality of a typical strongly interacting
subsystem, but the densities of the different subsystems are nonlinearly coupled.

In this paper our goal is to
\begin{itemize}
\item derive the mean field model,
\item validate it both by establishing exact properties and comparing to simulations of the full problem in low-dimensional examples.
\end{itemize}

The theoretical and computational results appear to us to be extremely promising.
Our theoretical results include mass conservation, energy conservation, and asymptotic
correctness in the limit of weak subsystem coupling (see Section~\ref{sec:mf}). The numerical tests we performed
on small systems show excellent qualitative agreement between mean field and exact model even when the coupling
is of order one,
at greatly reduced computational cost (see Section~\ref{sec:numerics}).

A fuller theoretical explanation of this good performance even for order one coupling would
be highly desirable, but lies beyond the scope of the present paper. 
One important aspect of our simulations appears to be a careful choice of subsystems which makes at least the {\it potential}
part of the Hamiltonian non-interacting. Remarkably, such a choice is always possible
in chain molecules with a standard force field (consisting
of nearest-neighbour bond terms, third neighbour angular terms,
and fourth neighbour torsion terms), by working in inner coordinates (i.e. bond lengths,
bond angles and torsion angles). In these coordinates, subsystem coupling only occurs through 
momentum exchange. But physical considerations as well as
previous work support the belief that fine details of momentum transfer do not play a decisive role in conformation dynamics:
\\
-- Molecules in solution are subject to relentless perturbations of momenta due to collisions with solvent molecules, 
yet conformation dynamics robustly takes place under these conditions.
\\
-- Many of the standard models in conformation dynamics involve randomization of momenta. In the Langevin equation,
this happens by addition of white noise to the momentum equation; in the approach by Sch\"utte in \cite{Sc99a}, 
one thinks of molecular conformations as subsets of \emph{configuration} space, 
considers only a \emph{spatial} transfer operator, 
and draws momenta at random from their statistical distributions in each computational evolution step. 

A more detailed analysis as well as applications of the mean field model to large systems are currently in 
progress and will appear elsewhere.

\section{Hamiltonian dynamics and Liouville equation}
In situations when quantum effects can be neglected and no bond-breaking or bond-formation
takes place, the dynamics of a molecule with $N$ atoms moving about in $\R^3$ can be described 
by a Hamiltonian of form  
\be \label{ham}
   H(q,p) = \frac{1}{2}p\cdot M(q)^{-1} p + V(q),
\ee
where $(q,p)\in\R^{2d}$, the mass matrix $M$ is a positive $d\times d$ Matrix, 
and $V\, : \, \R^d\to\R$ is a potential describing the atomic interactions.

In the case when all degrees of freedom are explicitly included and cartesian coordinates are used, 
we have $d=3N$ (where $N$ is the number of atoms), 
$q=(q_1,\ldots,q_N)\in\R^{3N}$, $p=(p_1,\ldots,p_N)$, and $M=\diag(m_i I_{3\times 3})$, 
where $q_i\in\R^3$, $p_i\in\R^3$, $m_i>0$ are the position, momentum, and mass
of the $i^{th}$ atom. In this paper we work with the more general form (\ref{ham}), in which the
kinetic energy is a quadratic form of $p$ depending on $q$. This form arises when inner coordinates
are used, which will play an important role below. For an $N$-atom chain molecule, the latter
consist of the $(N-1)$ nearest neighbour bondlengths $r_{ij}$, the $(N-2)$ bond angles $\theta_{ijk}$
between any three successive atoms, and the $(N-3)$ torsion angles $\phi_{ijkl}$ 
between any four successive atoms. In order to accurately model conformation changes, $V$
will have to contain at least nearest-neighbour bond terms $V_{ij}(r_{ij})$, third neighbour angular terms
$V_{ijk}(\theta_{ijk})$, and fourth neighbour torsion terms $V_{ijk\ell}(\phi_{ijk\ell})$. 
In practice the potentials could either come from
a suitable semi-empirical molecular force field model or from ab-initio computations.

The Hamiltonian dynamics takes the form
\begin{subequations} \label{MD}
\begin{align}
   \dot{q} & =  \phantom{-}\frac{\partial H}{\partial p}(q,p) = M(q)^{-1}p,\\ 
   \dot{p} & =  -\frac{\partial H}{\partial q}(q,p) = - \frac{\partial}{\partial q}\left(\frac{1}{2}p\cdot M(q)^{-1}p\right) 
 - \nabla V(q).
 \end{align}
\end{subequations}
It will be convenient to denote the phase space coordinates by $z=(q,p)\in\R^{2d}$ and the Hamiltonian
vector field by 
\be \label{hamvf}
f:= \left(\begin{array}{r} 
	\frac{\partial H}{\partial p} \\ 
	-\frac{\partial H}{\partial q} 
  \end{array}\right), 
\ee
so that (\ref{MD}) becomes 
\be \label{MDshort}
   \dot{z} = f(z).
\ee 
The Liouville equation associated to (\ref{MD}) describes the transport of a passive scalar $u$
by the flow associated to eq.~(\ref{MD}): 
\be \label{Liouville}
   \partial_t u + f \cdot \nabla_z u = 0,
\ee
where $u=u(z,t)$, $u\, : \, \R^{2d}\times\R\to\R$. 
Because the Hamiltonian vector field $f$ is divergence-free,
eq.~(\ref{Liouville}) can be written equivalently in the form 
\be \label{continuity}
   \partial_t u + div_z(u\, f) = 0.
\ee
Since eq.~(\ref{continuity}) preserves both positivity of $u$ and the total integral $\int_{\R^{2d}}u(z,t)\, dz$, 
it defines an evolution on the space of probability densities 
\[
\left\{u\in L^1(\R^{2d}) \, | \, u\ge 0, \int u = 1\right\}.
\]
Physically it can be interpreted as an evolution equation for ``ensembles'' of initial data under (\ref{MD}). The evolution law (\ref{MD}) for ``sharp'' initial data can be recovered
as a special case: $z(t)$ is a solution to (\ref{MD}) if and only if the delta function transported along this solution, $u(z,t)=\delta_{z(t)}$, is a solution to (\ref{continuity}).

Finally we note that eq.~(\ref{continuity}) preserves the (expected value of) energy, 
$$
   E(t) := \int H(z) u(z,t)\, dz.
$$
This is because by an integration by parts
$$
   \frac{d}{dt}E(t) = \int H(z)\Bigl( -\, div(u(z,t)f(z))\Bigr)\, dz 
   = \int \nabla H(z) \cdot f(z) u(z,t)\, dz
$$
and the inner product $\nabla H(z)\cdot f(z)$ vanishes for all $z$, due to (\ref{hamvf}). 

\section{Molecular conformations and almost invariant sets}\label{sec:conf}

A conformation of a molecule -- as we understand it \cite{DeDeJu99a,Sc99a,ScHuDe01a}	-- is given by an \emph{almost invariant} (or \emph{metastable}) subset of configuration space.  

Roughly speaking, an almost invariant set of a discrete dynamical system $S:X\to X$ is a subset $A\subset X$ such that the \emph{invariance ratio}
\[
\rho_\mu(A) = \frac{\mu(S^{-1}(A)\cap A)}{\mu(A)}
\]
is close to $1$, cf.\ \cite{DeJu99a}. Here $\mu$ denotes a suitable probability  measure, typically Lebesgue measure or some $S$-invariant measure.  In our case, $S$ will be the time-$T$-map $\Phi^T:\R^{2d}\to\R^{2d}$ of the Hamiltonian system (\ref{MD}) (More generally, instead of a deterministic map, one can consider a stochastic process, cf.\ \cite{DeJu99a}.) 

Almost invariant sets can be determined via the computation of the eigenfunctions at (real) eigenvalues close to one of a certain \emph{transfer operator} associated to $S$ \cite{DeJu99a,Sc99a}.  For example, the \emph{Frobenius-Perron} operator 
\[
P\mu(A) = \mu(S^{-1}(A)), \quad A\text{ measurable},
\]
describes the evolution of measures on phase space.  This is a linear operator on the space $\cM_\C$ of bounded complex valued measures on $X$. By definition, $\|P\mu\| \leq \|\mu\|$ and thus the spectrum of $P$ is confined to the unit circle.  Eigenmeasures $\mu(A)=\mu(S^{-1}(A))$ at the eigenvalue $1$ are \emph{invariant measures}.  If $0\neq\mu\in\cM_\C$ is an eigenmeasure of $P$ at the eigenvalue $\lambda$, then
\[
\lambda \mu(X) = P\mu(X) = \mu(S^{-1}(X)) = \mu(X)
\]
and thus $\mu(X)=0$ for $\lambda\neq 1$. In particular, if $\lambda< 1$ and $\mu$ are real, then there are two positive real measures $\mu^+, \mu^-$ such that $\mu=\mu^+ - \mu^-$ (\emph{Hahn-Jordan decomposition}).  If $\mu$ is normalized such that $|\mu|=\mu^+ +\mu^-$ is a probability measure then \cite{DeJu99a}
\[
\rho_{|\mu|}(A^+) + \rho_{|\mu|}(A^-) = \lambda + 1
\]
where $A^+=\supp(\mu^+)$ and $A^-=\supp(\mu^-)$.  Consequently, if $\lambda < 1$ is close to $1$ then both $\rho_{|\mu|}(A^+)$ and $\rho_{|\mu|}(A^-)$ are close to $1$ and thus almost invariant. 
\paragraph{Transfer operators.}

For measures which are absolutely continuous one can equivalently consider $P$ on $L^p(X,\C)$. Since the flow $\Phi^T$ of (\ref{MD}) is a volume-preserving diffeomorphism, this operator takes a particularly simple form in our case\footnote{we write $P^T$ in order to stress the dependence of the operator on the integration time $T$}:
\begin{equation}\label{eq:Liou_time-T}
P^T u = u\circ\Phi^{-T},
\end{equation}
which is the time-$T$-map of the Liouville equation (\ref{Liouville}). Note that for an arbitrary 
function $g:\R\to [0,\infty)$ of the Hamiltonian, the function $u(z)=g(H(z))$ satisfies $\nabla_z u(z)=g'(H(z))\nabla_z H(z)$.  Thus $f\cdot\nabla_z u = 0$ and $u$, normalized s.t.\ $\int u(z)\, dz = 1$, is an invariant density.  Of particular interest is the \emph{canonical density}
\be\label{eq:canonical density}
h(z)=C\exp(-\beta H(z)),
\ee
$C=\int \exp(-\beta H(z))\, dz$, where $\beta=1/(kT)$ and $k$ is Boltzmann's constant.  This density describes the distribution of a (constant) large number of molecules at temperature $T$ and of constant volume.  Note that we can write
\[
h(z)=h(q,p)=C_p(q)\exp\left(-\frac{\beta}{2}p\cdot M^{-1}(q)p\right)C_q\exp\left(-\beta V(q)\right)=:h_p(q,p)h_q(q),
\]
where $C_p(q)$ and $C_q$ are chosen such that $\int h_p(q,p)dp = 1$ for each $q$, and $\int h_q = 1$.
\paragraph{Spatial transfer operator.}

As mentioned, molecular conformations should be thought of as almost invariant subsets of configuration space.  Sch\"utte \cite{Sc99a} introduced a corresponding spatial transfer operator by averaging (\ref{eq:Liou_time-T}) over the momenta:  Let $h\in L^1(\R^{2d})$ be an invariant density of (\ref{eq:Liou_time-T}) with $h(q,p)=h(q,-p)$, let $\bar h(q)=\int h(q,p)\, dp$ and consider the operator  
\begin{equation}\label{eq:P_spatial}
S^Tu(q) = \frac{1}{\bar h(q)}\int u\left(\pi_q\Phi^{-T}(q,p)\right)\, h(q,p)\, dp,
\end{equation}
where $\pi_q(q,p)=q$ is the canonical projection onto configuration space.  Sch\"utte \cite{Sc99a} showed that under suitable conditions, the spatial transfer operator is self-adjoint and quasi-compact on an appropriate weighted $L^2$ space.

\paragraph{Transition probabilities.}

A key quantity of interest are transition probabilities  from one region of space into another.  The transition probability from a region $B_i\subset\R^{2d}$ into another region $B_j\subset\R^{2d}$ in phase space is given by the volume fraction of those initial data in
$B_i$ which end up in $B_j$ at time $T$,
\be \label{prob}
   p^{(T)}_{ij} = \frac{m(\Phi^{-T}(B_j)\cap B_i)}{m(B_i)} = \frac{1}{m(B_i)} \langle P^T\chi_{B_j},\chi_{B_i}\rangle
\ee
where $m$ denotes $2d$-dimensional volume.   Using the expression based on the transfer operator, similarly transition probabilities between subsets of configuration space can be defined.  Typically, these quantities are sought for large $T$, but only short time simulations are numerically feasible.  However, $p^{(T)}_{ij}$ can be approximated for large $T$ by repeated matrix-vector multiplications, requiring short time evaluations of the flow only, cf.\ \cite{DeJuKo05a}.

\section{Mean field approximation} \label{sec:mf}

The problem with eq.~(\ref{continuity}) as it stands is that it is amenable to a direct numerical treatment only
for very small systems,
due to the exponential scaling of the number of computational degrees of freedom with particle number.
If the phase space of each atom is approximated by a $K$-point grid, such that the solution to 
(\ref{continuity}) at time $t$ becomes a vector in $\R^K$, then the corresponding grid of the $N$-particle
system has $K^N$ gridpoints and the solution at time $t$ becomes a vector in $\R^{K^N}$. Our proposal
to address this problem is partially inspired by Hartree-Fock- and density functional theory methods
in quantum chemistry, which allow to overcome a related complexity problem for the $N$-particle Schr\"odinger
equation. 
\paragraph{Partitioning into subsystems.} Starting point is an, for the moment arbitrary, partition of
phase space coordinates $z=(q,p)$ into subsystem coordinates: 
$$
   z=(z_1,\ldots,z_N)\in\R^{2d}, \;\;\; z_i=(q_i,p_i)\in\R^{2d_i}, \;\;\; \sum_{i=1}^Nd_i=d,
$$
where $p_i$ is the vector of momentum coordinates corresponding to the position
coordinates $q_i$. Let $f_i=\left(\frac{\partial H}{\partial p_i}, \, -\frac{\partial H}{\partial q_i}\right)$. Then eq. (\ref{MD}) can be
re-written as
\be \label{MDshort'}
   \dot{z_i} = f_i(z), \;\;\; i=1,..,N.
\ee 

Given a phase space density $u(z_1,\ldots,z_N,t)$ which
depends on the phase space coordinates of all the subsystems, we introduce reduced densities for each
subsystem $i=1,\ldots,N$, as follows: 
\be \label{densities}
   u_i(z_i,t) := \int_{\R^{2(d-d_i)}}u(z,t)\, d\zihat, \\
\ee
where here and below $\zihat$ denotes the coordinates $(z_j)_{j\neq i}$. 
Note that, by the normalization of $u$, 
$$
    \int_{\R^{2d_i}}u_i(z_i,t)\, dz_i = 1,
$$
that is to say the probability that the $i^{th}$ particle has {\it some} position and momentum is one.
Next we calculate the exact time evolution of the subsystem densities $u_i$.  
By (\ref{continuity}) and (\ref{densities}), 
\begin{eqnarray*}
   \partial_tu_i(z_i,t) &=& -\int div_{z}(u(z,t)\, f(z)) \, d\widehat{z_i} 
                         = -\int div_{z_i} (u(z,t)\, f_i(z)) \, d\widehat{z_i} \\
                        &=& -\, div_{z_i} \int u(z,t)\, f_i(z)\, d\widehat{z_i}.
\end{eqnarray*}
Using (\ref{densities}) this can be rewritten in the form of a single-subsystem transport equation,
\be \label{exact}
   \partial_tu_i(z_i,t) = - div_{z_i} \Bigl(u_i(z_i,t) \tilde{f}_i(z_i,t)\Bigr),
\ee
with underlying vector field
\be \label{exact2}
    \tilde{f}_i(z_i,t) = \frac{\int u(z,t)\, f_i(z)\, d\widehat{z_i}}{\int u(z,t)\, d\zihat}.
\ee
Note, however, that (\ref{exact}), (\ref{exact2}) is not a closed system, since the vector
field $\tilde{f}_i$ depends on the full density $u$, not just the subsystem densities $u_j$ ($j=1,..,N$).

\paragraph{Mean field model.} To close the system (\ref{exact}), (\ref{exact2}), we now make the statistical 
independence ansatz 
\be \label{ansatz}
   u(z_1,\ldots,z_N,t) = u_1(z_1,t)\cdots u_N(z_N,t).
\ee
The vector field $\tilde{f}_i(z_i,t)$ in (\ref{exact}) then simplifies to 
\be \label{F1}
     \tilde{f}_i(z_i,t) = \int_{\R^{2(d-d_i)}} f_i(z) \prod_{j\neq i} u_j(z_j,t)\, d\zihat =: f_i^\text{mf}\Bigl[\widehat u_i\Bigr](z_i,t). 
\ee
We call the system of eqs.~(\ref{exact}), $i=1,\ldots,N$, with $\tilde{f}_i$ given by (\ref{F1})
the {\it mean field approximation} to
the Liouville equation. Note that it is a system of $N$ coupled nonlinear partial integrodifferential
equations on the lower-dimensional subsystem phase spaces $\R^{2d_i}$, whereas the original Liouville equation was a 
linear partial differential equation on $\R^{2d}$, $d=\sum_i(2d_i)$. Physically, each subsystem can be pictured, at time $t$, as experiencing the force of
  the ensemble of the other subsystems in their ``typical'' states at time $t$.

We record some basic properties of the mean field approximation.
\begin{enumerate}
\item The total densities $\int u_i(z_i,t)\, dz_i$ are conserved. This is immediate from the conservation law form
$\partial_tu_i+div\, (u_i f_i^\text{mf})=0$. Thus we may assume $\int u_i(z_i,t)\, dz_i=1$ for all $i$ and $t$, and continue to 
interpret the $u_i$ as probability densities. 

\item For non-interacting subsystems, i.e. 
$$
   H(z)=\sum_{i=1}^N\Bigl( \frac12 p_i\cdot M_i(q_i)^{-1}p_i + V_i(q_i)\Bigr),
$$ 
it is exact, that is to say if the $u_i(z_i,t)$ evolve via (\ref{exact}), (\ref{F1}), then
the product $u_1(z_1,t)\cdots u_N(z_N,t)$ solves the original Liouville equation (\ref{continuity}). This follows
from the fact that in this case, 
$$
   \frac{\partial H}{\partial p_i} = M_i(q_i)^{-1}p_i, \;\;\; -\frac{\partial H}{\partial q_i} = 
   - \frac{\partial}{\partial q_i}\mbox{$\frac12$}p_i\cdot M_i(q_i)^{-1}p_i 
   - \nabla V_i(q_i)
$$
and so $f_i(z)$ depends only on $z_i$, as a consequence of which the exact vectorfield  $f_i$, the reduced
vectorfield $\tilde{f}_i$ in (\ref{exact2}), and the mean field vectorfield $f_i^\text{mf}$ (RHS of (\ref{F1})) all coincide,
regardless of the mean field ansatz (\ref{ansatz}). 

\item For given $u_j$, $j\neq i$, the transport equation for $u_i$ has the form of a Liouville equation coming from
an associated time-dependent subsystem Hamiltonian, 
\be \label{ea:td_ham}
    H_i^\text{mf}(q_j,p_j,t) = \int H(q,p) \prod_{j\neq i} u_j(q_j,p_j,t) \, d\zihat,
\ee
that is to say
\be \label{hamform}
   f_i^\text{mf}(q_i,p_i,t) = \left(\begin{array}{c} \frac{\partial}{\partial p_i} H_i^\text{mf}(p_i,q_i,t) \\
                    -\frac{\partial}{\partial q_i} H_i^\text{mf}(p_i,q_i,t) \end{array}\right).
\ee
In particular, $f_i^\text{mf}$ is divergence-free. 
Note that time-dependence of the effective subsystem Hamiltonian enters only through time-dependence of
the $u_j$, $j\neq i$.    

\item The total energy 
$$
   E(t) := \int H(z)\, u_1(z_1,t)\cdots u_N(z_N,t)\, dz_1\cdots dz_N
$$
is conserved. To see this, calculate using the evolution equation for the $u_i$ and an integration by parts
\begin{eqnarray*}
  \frac{d}{dt} E(t) &=& \int H(z) \left[\sum_{i=1}^N \frac{\partial u_i}{\partial t}(z_i,t)\prod_{\ell\neq i}
                        u_\ell(z_\ell,t)\right]dz_1\cdots dz_N \\
  &=& \int H(z) \left[\sum_{i=1}^N -\, div_{z_i}\left(u_i(z_i,t) f_i^{\text{mf}}(z_i,t)\right)
                        \prod_{\ell\neq i} u_\ell(z_\ell,t)\right] dz_1\cdots dz_N \\
  &=& \sum_{i=1}^N \int \left[ \int \nabla_{z_i}H(z) \prod_{\ell\neq i} u_\ell(z_\ell,t)\, d\widehat{z_i}\right]
      \cdot u_i(z_i,t)\,f_i^\text{mf}(z_i,t)\, dz_i.
\end{eqnarray*}
But the term in square brackets equals $\nabla_{z_i}H_i^{\text{mf}}(z_i,t)$, and hence its
dot product with $f_i^{\text{mf}}(z_i,t)$ vanishes for all $z_i$ and $t$, on account of (\ref{hamform}); consequently
$\frac{d}{dt} E(t)=0$. 
\end{enumerate}

Property 2.\ contains useful information regarding how the, up to now arbitrary, partitioning into
subsystems should be chosen in practice. In order to maximize agreement with the full Liouville equation
(\ref{Liouville}), the subsystems should be only weakly coupled. In the case of an $N$-atom chain,
this suggests to work with subsystems defined by inner, not cartesian, coordinates (as is done in the simulation of n-butane
in Section~\ref{subsec:butane} below). Namely, in inner coordinates,
at least the {\it potential energy} decouples completely for
standard potentials containing nearest-neighbour bond terms, third neighbour angular terms and fourth
neighbour torsion terms: $V((r_{ij}), \,(\theta_{ijk})_{ijk},\,(\phi_{ijk\ell})_{ijk\ell})=\sum V_{ij}(r_{ij}) + \sum
V_{ijk}(\phi_{ijk}) + \sum V_{ijk\ell}(\phi_{ijk\ell})$.

A deeper, and perhaps surprising, theoretical property of the mean field model which goes beyond property 2.\ concerns weakly coupled
subsystems. Consider a Hamiltonian of form $H(z)=H_0(z)+\epsilon H_{{\rm int}}(z)$, where
$H_0$ is a non-interacting Hamiltonian of the form given in 2., and $\epsilon$ is a coupling constant.
It can then be shown that the exact subsystem densities
$u_i$ obtained from (\ref{Liouville}), (\ref{densities}) and the mean field densities obtained by solving
(\ref{exact}), (\ref{F1}) differ, up to any fixed time $T$, only by $O(\epsilon^2)$, not the naively expected $O(\epsilon)$. This means
that the effect of coupling between subsystems is captured correctly to leading order (in the coupling constant)
by the mean field approximation. For a proof of this fact see our companion paper \cite{FrJuKo09a}.

\section{The mean field transfer operator}

The mean field approximation to the Liouville equation introduced above gives rise in a natural way to a mean field approximation of the transfer operator
$$
   P^T u^0 = u(\cdot, T), 
$$
where $u$ is the solution to the Liouville equation (\ref{Liouville}) with initial condition
$u(z,0)=u^0(z)$.  We define the \emph{mean field transfer operator} as
\[
P^T_\text{mf}(\bu^0) = \bu(\cdot,T),
\]
where $\bu(z,t)=(u_1(z_1,t),\ldots,u_N(z_N,t))$ is the solution to the mean field approximation (\ref{exact}), (\ref{F1}), $i=1,\ldots,N$, of the Liouville equation with initial data $\bu(z,0)=\bu^0(z)=(u^0_1(z_1),\ldots,u^0_N(z_N))$.  Note that the components of $P^T_\text{mf}$ are multilinear, i.e.\ for fixed $\widehat \bu_i$, the map
\[
P^T_{\text{mf}}(\widehat \bu_i)u_i := \left[P^T_\text{mf}(\bu)\right]_i 
\]
is linear.
 
\paragraph{The mean field spatial transfer operator.}

Consider an ensemble of molecules whose distribution in phase space is given by an invariant density $h(q,p)$. We would like to describe distribution changes in the position space only (cf.~Section~\ref{sec:conf}).  Following \cite{Sc99a} we define the \textit{spatial transfer operator} 
\begin{equation}
	S^Tw(q) = \int P^T\left(w(q)\bar{h}(q,p)\right)\intd p,
	\label{eq:spatial_transop}
\end{equation}
where $\bar{h}$ is the conditional density of $p$ for a given $q$, i.e.
\[ \bar{h}(q,p) = \frac{h(q,p)}{\int h(q,p)\, dp}. \]
Now we define the spatial transfer operator corresponding to the mean field system. The distribution of the $i$-th subsystem is given by
\[ h_i(z_i) = \int h(z)\intd \widehat z_i. \]
The conditional density of the $p_i$ for given $q_i$ is
\[ \bar{h}_i(q_i,p_i) = \frac{h_i(q_i, p_i)}{\int h_i(q_i,p_i)\, dp_i}. \]
We therefore define the \emph{mean field spatial transfer operator} as
\begin{equation}
	S_\text{mf}^T (\bw^0) = \bw(\cdot, T),
\end{equation}
where $\bw^0(q)=(w^0_1(q_1),\ldots,w^0_N(q_N))$, $\bw(q,t)=(w_1(q_1,t),\ldots,w_N(q_N,t))$ and
\begin{equation}	
		w_i(q_i,T) = \int \left[P_\text{mf}^T(\bu^0)\right]_i(q_i,p_i)\intd p_i
	\label{eq:spatial_mf_1}
\end{equation}
with $u^0_i(q_i,p_i) = w^0_i(q_i)\bar h(q_i,p_i)$, $i=1,\ldots,N$.  Again, the components of $S^T_\text{mf}$ are multilinear, i.e.\
for fixed $\widehat \bw_i$, the map
\be\label{eq:MFST_ml}
S^T_{\text{mf}}(\widehat \bw_i)w_i := \left[S^T_\text{mf}(\bw)\right]_i
\ee
is linear.

\paragraph{Mean field eigenfunction approximation.}

Our goal is to compute eigenfunctions of the (full) spatial transfer operator (\ref{eq:spatial_transop}). As approximations, we are going to use products of eigenfunctions (at the leading eigenvalues) of the linear component maps (\ref{eq:MFST_ml}).  In computing these, we fix $\widehat \bw_i$ to the invariant density of $S^T_\text{mf}$.  This approach is motivated by the following observation for non-interacting subsystems: Let $S_1,S_2:X\to X$ be two maps and $P_1,P_2:L^1\to L^1$ the associated Frobenius-Perron operators.  For the product map $S=S_1\otimes S_2:X^2\to X^2$, $S(x_1,x_2)=(S_1(x_1),S_2(x_2))$, the Frobenius-Perron operator is given by
\[
(Ph)(x_1,x_2) = (P_1h_1)(x_1)\cdot (P_2h_2)(x_2),
\]
if $h(x_1,x_2) = h_1(x_1)h_2(x_2)$. Now let $P_1h_1 = \lambda_1h_1$ and $P_2h_2 = \lambda_2h_2$ for some eigenvalues $\lambda_1,\lambda_2\in\C$, then $\lambda_1\lambda_2$ is an eigenvalue of $P$ with eigenfunction $h_1h_2$.

In order to compute the invariant density, we resort to an iterative procedure which is inspired by the famous Roothaan algorithm from quantum chemistry: For a given initial guess $\bw^0$ we iteratively compute a sequence $\bw^k$, $k=0,1,\ldots$ of approximate invariant densities of $S^T_\text{mf}$ by computing the fixed point of each linear component map $S^T_\text{mf}(\widehat \bw_i^k)$, i.e.\ by computing $\bw^{k+1}=(w^{k+1}_1,\ldots,w^{k+1}_N)$ such that 
\[
S^T_\text{mf}(\widehat \bw^k_i)w^{k+1}_i = w^{k+1}_i, \quad i=1,\ldots,N.
\]

\section{Numerics and examples}\label{sec:numerics}

\subsection{Ulam's method}

The (spatial) transfer operator is typically considered as an operator on a (suitably weighted) $L^p$ ($p=1$ or $2$) space.  As such, it is not directly amenable to numerical computations. In \cite{Ul60a}, Ulam proposed the following discretization: Let $\mathcal X_n = \{X_1,\ldots,X_n\}$ be a disjoint partition of (phase or configuration) space, $V_n := \text{span}\{\chi_{1},\ldots,\chi_{n}\}$, where $\chi_{i}$ denotes the characteristic function of $X_i$ and $Q_n$ the projection from $L^p$ onto $V_n$ defined by
\[ 
Q_{n}f:=\sum_{i=1}^{n}c_i\chi_{i}\qquad\text{ with } \qquad c_i:=\frac{1}{m(I_i)}\int_{I_i}f\intd m. 
\]
The \emph{discretized (spatial) transfer operator} $S^T:V_n\to V_n$ is defined to be
\[ 
S_n^T = Q_nS^T. 
\]
Note that $S_n^T$ can be represented by a stochastic matrix, where the matrix entries $(S_n^T)_{ij}$ are the transition rates between the sets $X_j$ and $X_i$. In other words, $(S_n^T)_{ij}$ is the probability, that $\pi_q\Phi^T(q,p)\in X_i$, if $q\in X_j$ is sampled according to a uniform distribution and then $p$ according to $\bar{h}(q,\cdot)$. 
This yields the numerical computation of the transition matrix by a Monte-Carlo method:
\[ (S_n^T)_{ij} \approx (\widetilde{S}_{n}^T)_{ij} :=  \frac{1}{K}\sum_{k=1}^K\chi_i\left(\pi_q\Phi^T\left(p^{(k)},q^{(k)}\right)\right), \]
where the points $q^{(1)},\ldots,q^{(k)}$ are chosen i.i.d. from $X_j$ according to a uniform distribution, and the correspondig $p^{(k)}$ according to $\bar{h}(q^{(k)},\cdot)$. Observe that the Markov-structure of the transfer operator is preserved: $S^T_\text{mf}$ and $Q_n$ preserve integral and positivity, and the Monte-Carlo approximation of the discretized transfer operator is a stochastic matrix, too.

Note that we have to evaluate the flow map several times for each partition element $X_i$, but that these are \textit{short} time simulations only. The time $T$ merely has to be large enough that motion can be observed. In particular, $T$ may be much smaller than characteristic times, e.g. for conformational changes; they can be several orders of magnitude apart. Thus, the problem of evaluating the flow at time $T$ is well conditioned, and we may use low order explicit schemes for the time integration.

\subsection{Complexity}

Let us first investigate the costs of setting up the discretized transfer operator for an arbitrary subsystem.
Using Ulam's method, we need to perform the following steps for each partition element $X_j$:
\begin{itemize}
	\item sample $q^{(k)}\in X_j$ and correspondingly $p^{(k)}$,
	\item integrate the mean field system for time $T$ and initial data $(p^{(k)},q^{(k)})$,
	\item project the endpoint onto position space and find the partition element $X_i$ it is contained in.
\end{itemize}
Using the canonical density for the invariant density $h$, there is an explicit representation for the momentum distribution $\bar{h}(q,\cdot)$ which can be sufficiently well approximated by a linear combination of Gaussians. The numerical time integration of the initial points requires several evaluations of the mean field vector field (\ref{hamform}).  This in turn requires the numerical evaluation of a $2(d-d_i)$ dimensional integral. The intergal w.r.t. the $\hat{p}_i$ can be handled analytically and by an apriori computation which is \textit{independent} of the $w_i$, $p_i$ and $q_i$.  Naively, this leaves us with a $d-d_i$ dimensional integral. However, note that in the case of non-interacting subsystems, i.e.\ $f_i(z)=f_i(z_i)$, $f_i$ can be pulled out and the integral reduces to 1.  For systems with small subsystems, i.e.\ $d_i\leq \bar d$ and $\bar d$ small, and in which only a fixed and small number of neighboring subsystems interact, the dimensionality of the integral is $\sum_{j\sim i} d_j = \mathcal{O}(\bar d)$, where $j\sim i$ means all subsystems $j$ which interact with subsystem $i$.

The solution of the resulting eigenvalue problems is simple compared with the assembling of the discretized mean field transfer operator, particularly since we are only interested in the dominant part of the spectrum. Arnoldi type iteration methods can be used.

\subsection{Example: a simple 2d-system}

\begin{figure}[htbp] 
   \centering
   \includegraphics[width=0.47\textwidth]{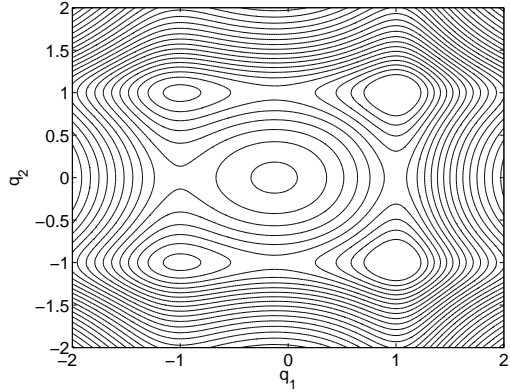} 
   \caption{Potential (\ref{V4}) of the simple 2d-system.}
   \label{fig:V4}
\end{figure}

For $(q_1,q_2)\in \R^2$ consider the potential
\begin{eqnarray}\label{V4}
V(q_1,q_2) &=& \left(\frac{3}{2} q_1^4+\frac{1}{4} q_1^3 -3 q_1^2
-\frac{3}{4}q_1+3\right)\cdot \left(2 q_2^4 -4q_2^2+\alpha\right)\\
&=& V_1(q_1)\cdot V_2(q_2)
\end{eqnarray}
with $\alpha=3$, cf.~Figure~\ref{fig:V4}.

\begin{figure}[tbp] 
   \centering
   \includegraphics[width=0.32\textwidth]{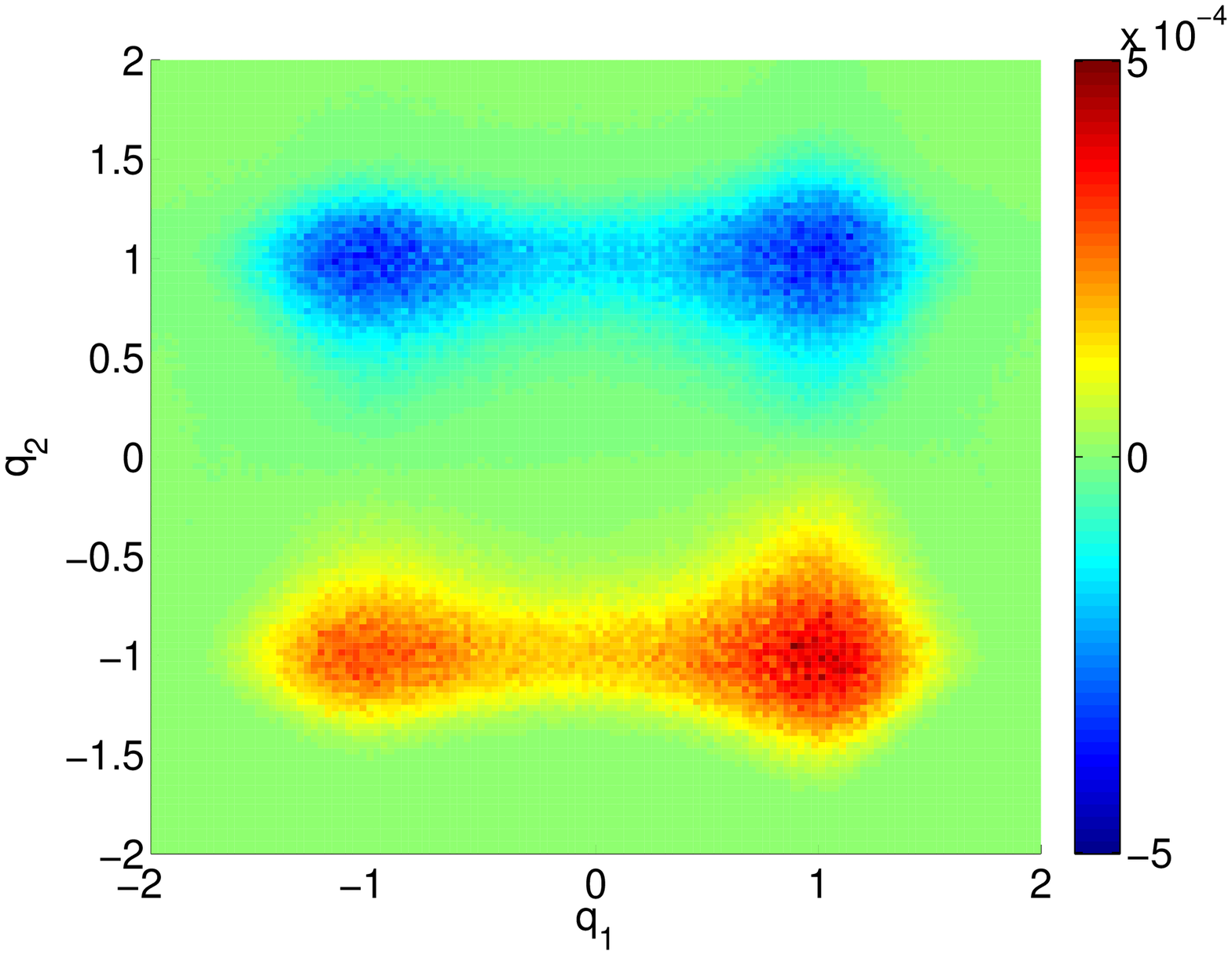} 
   \includegraphics[width=0.32\textwidth]{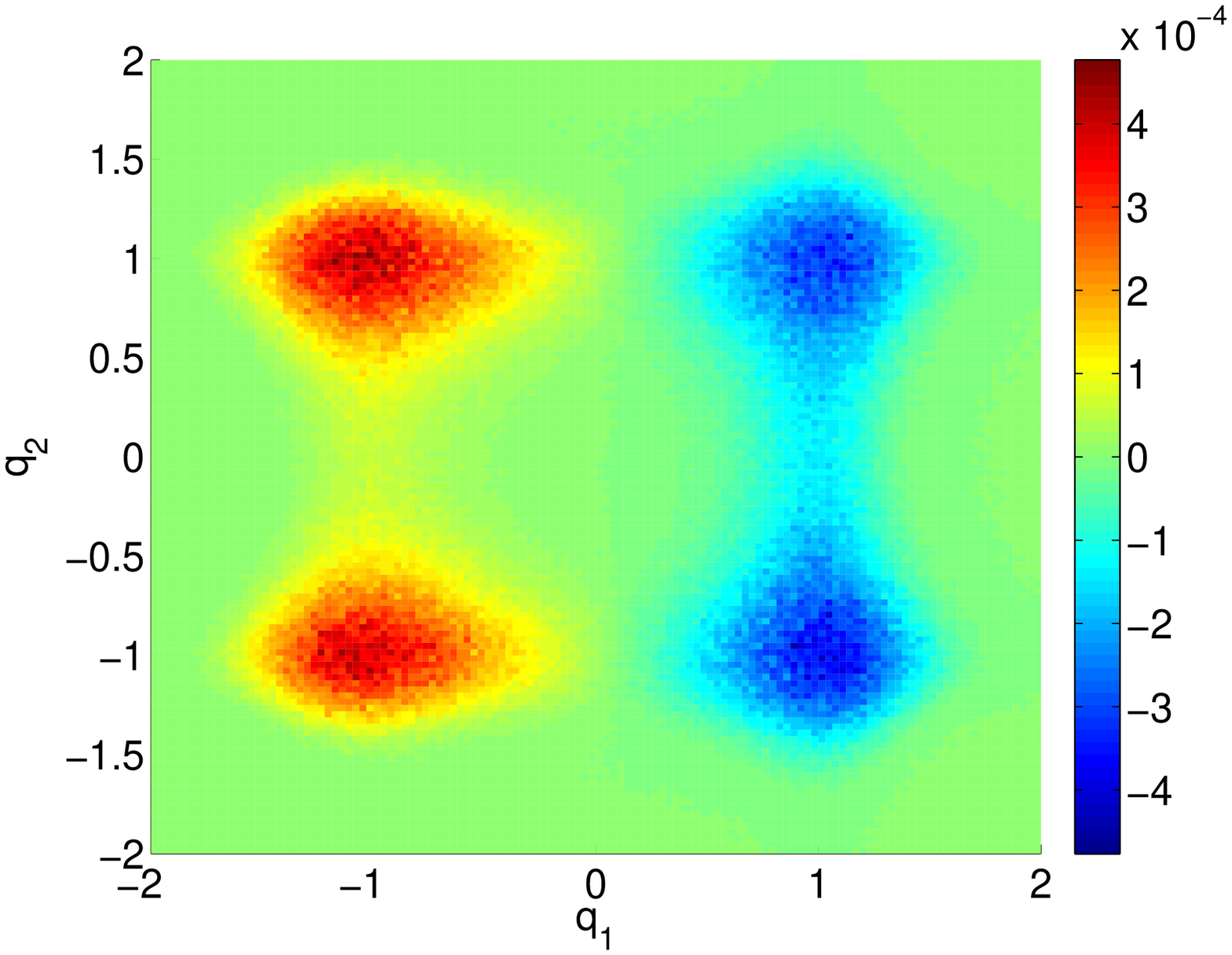} 
   \includegraphics[width=0.32\textwidth]{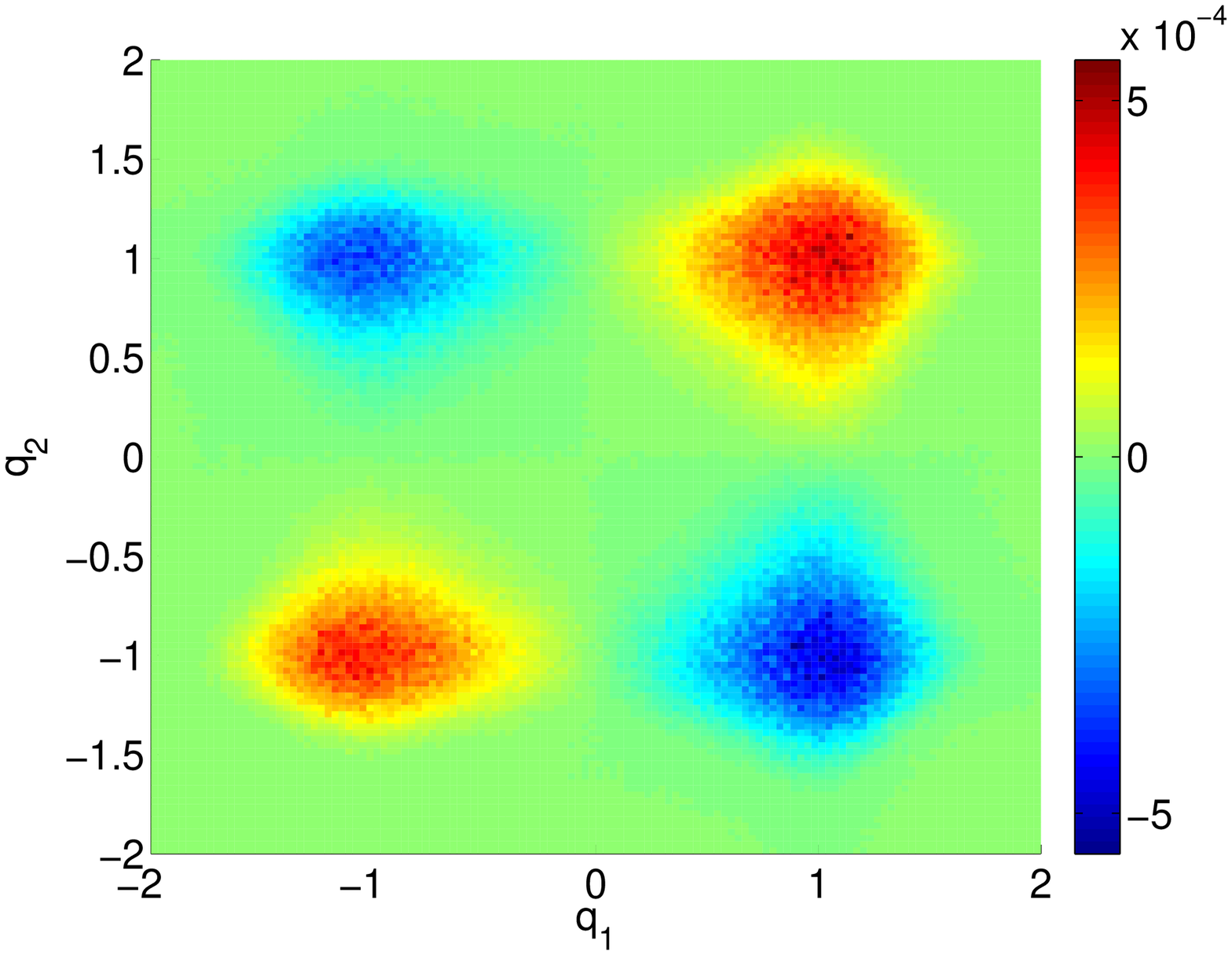} 
   \caption{Eigenfunctions at the second, third and fourth largest eigenvalue of the full spatial transfer operator (from left to right).}
   \label{fig:ex1_full}
\end{figure}

Figure~\ref{fig:ex1_full} shows the eigenfunctions at the leading eigenvalues of the full spatial transfer operator.

\begin{figure}[tbp] 
   \centering
   \includegraphics[width=0.37\textwidth]{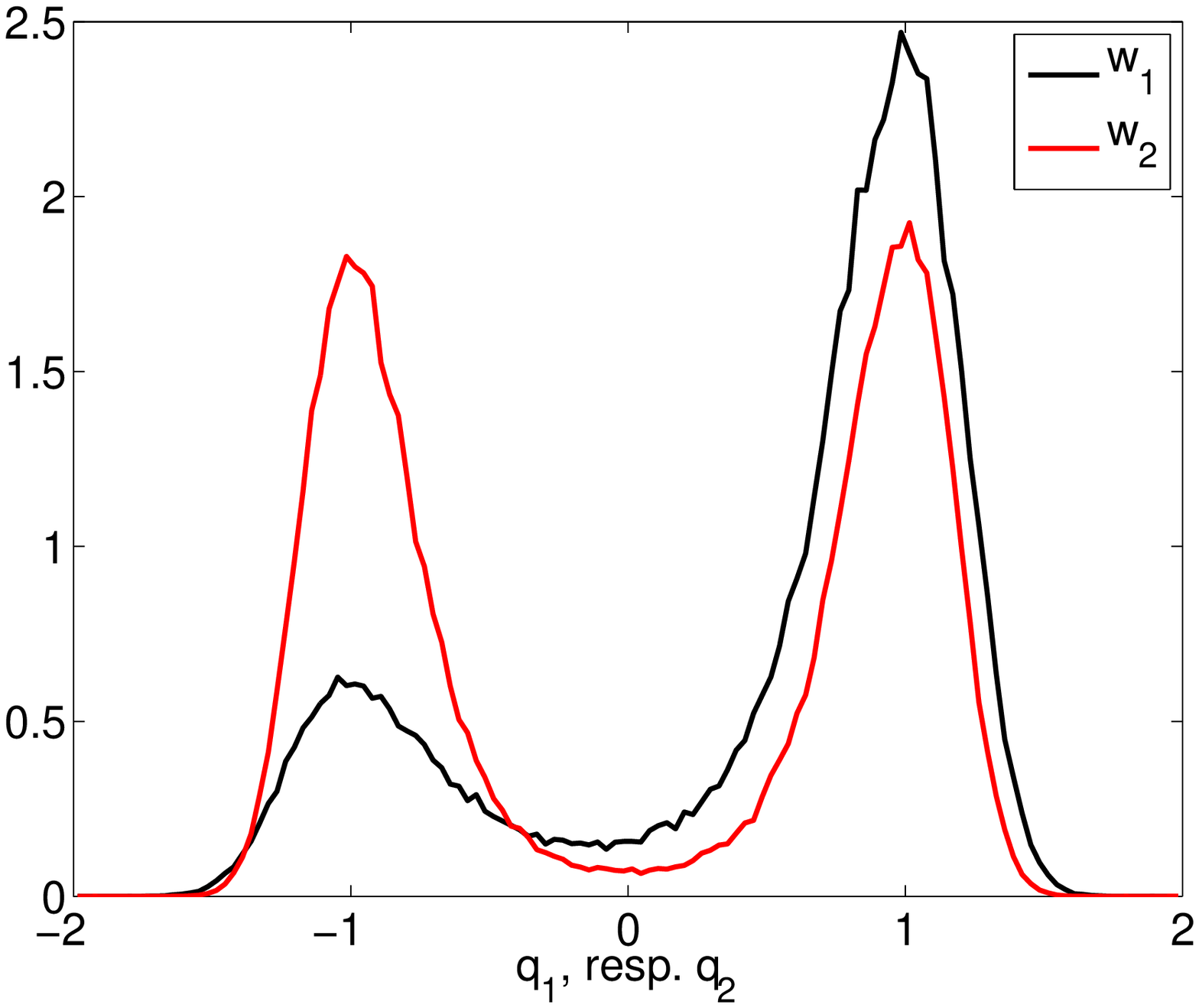}\qquad 
   \includegraphics[width=0.37\textwidth]{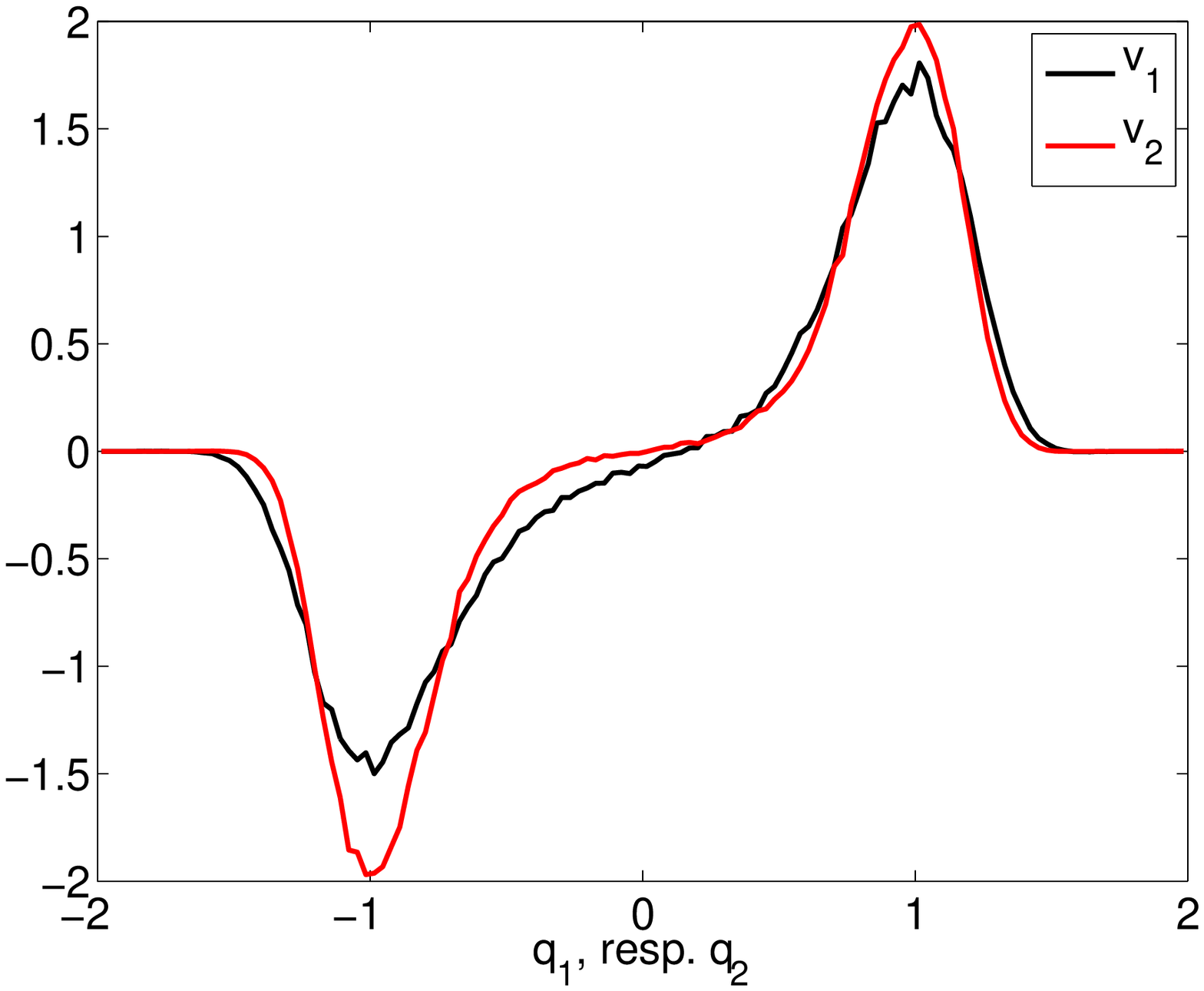} 
   \caption{Invariant density (left) and eigenfunctions at the second eigenvalue of the two components of the mean field spatial transfer operator.}
   \label{fig:ex1_factors}
\end{figure}

\begin{figure}[htbp] 
   \centering
   \includegraphics[width=0.32\textwidth]{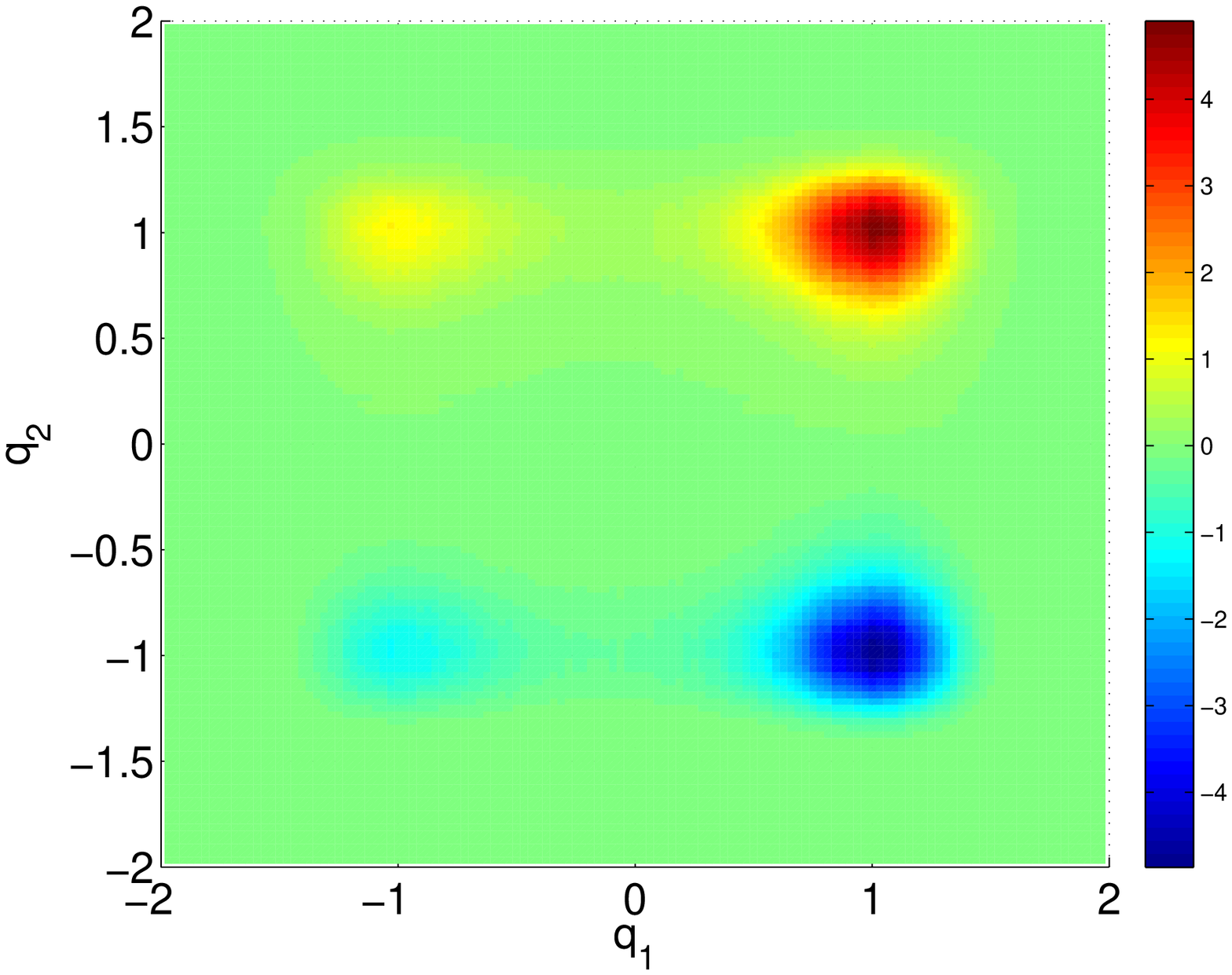} 
   \includegraphics[width=0.32\textwidth]{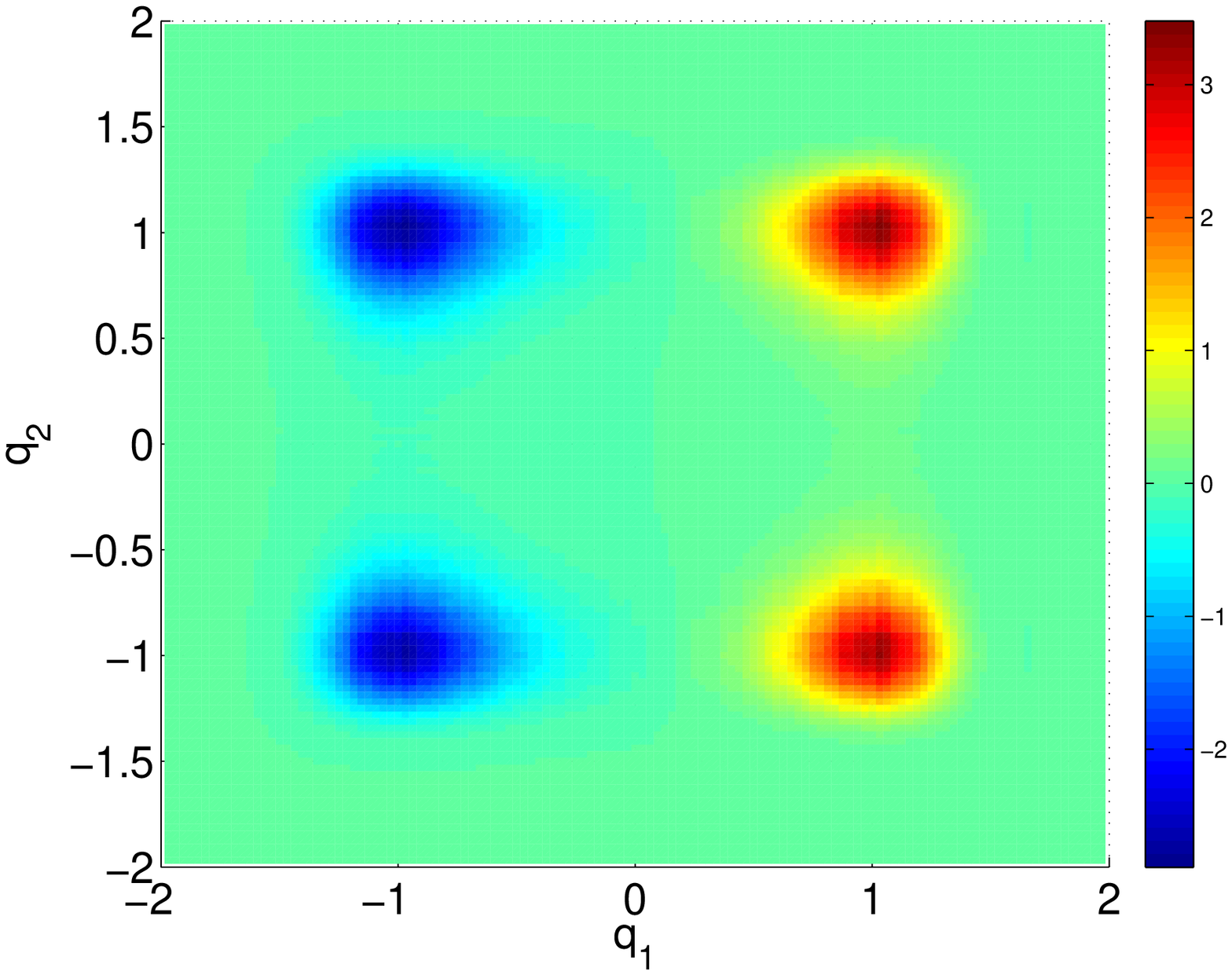} 
   \includegraphics[width=0.32\textwidth]{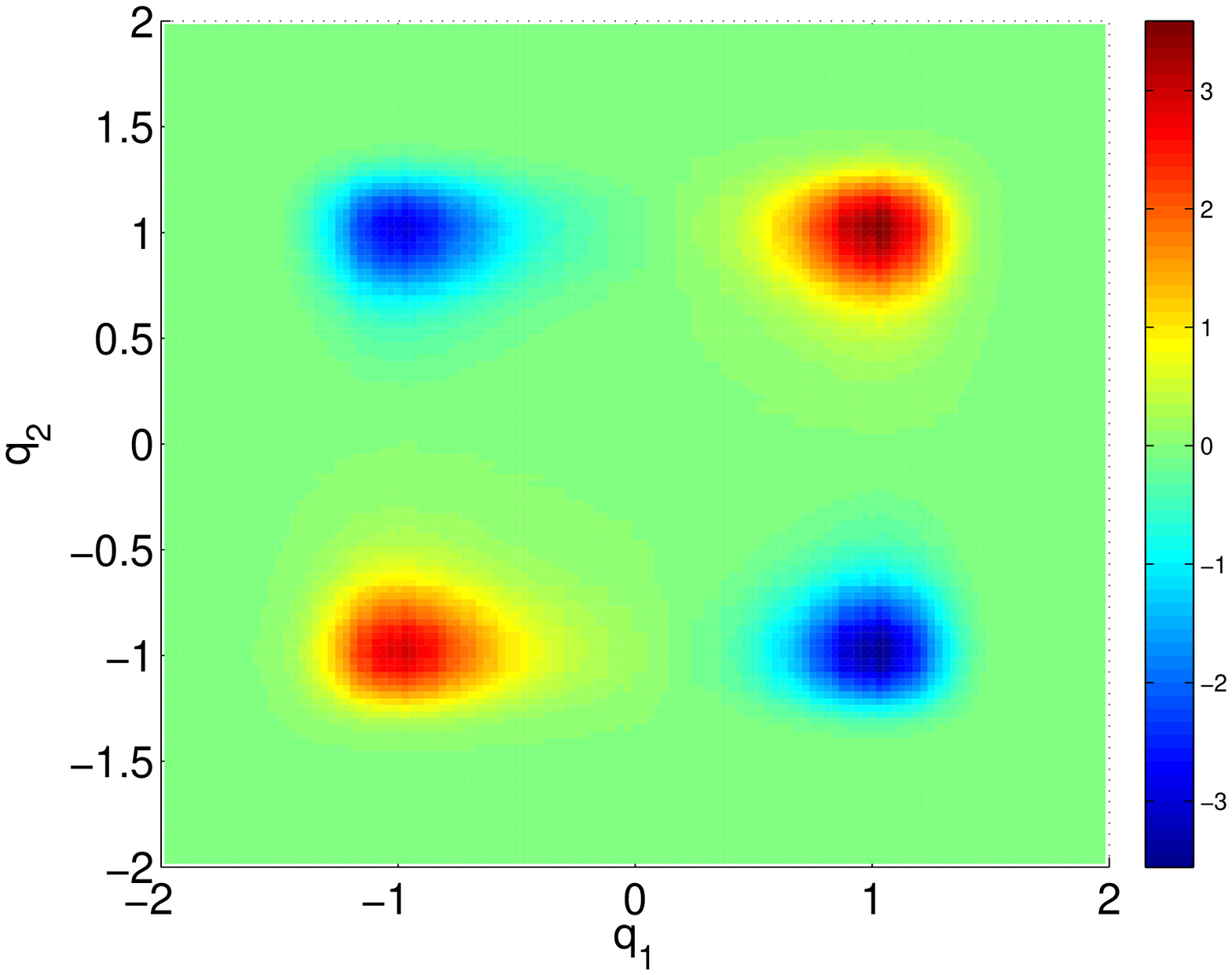} 
   \caption{Mean field approximations to the eigenfunctions from Fig.~\ref{fig:ex1_full}.}
   \label{fig:ex1_roothaan}
\end{figure}

The mean field approximation to the Liouville equation in this case reads explicitly
\[
   \partial_tu_i(z_i,t)=  \left(\begin{array}{c} -m_i^{-1}p_i \\
                    \displaystyle\nabla_{q_i}V_i(q_i)\int V_j(q_j)\, u_j(z_j,t)\, dz_j \end{array}\right) \cdot \nabla_{z_i} u_i(z_i,t), \quad i=1,2,\quad j\neq i.
\]
Figure~\ref{fig:ex1_factors} shows the invariant density $\bw=(w_1,w_2)$ (after ten iterations of the Roothaan type iteration) as well as the eigenfunctions $v_1$ and $v_2$ at the second eigenvalue of the two linear component maps $S^T_\text{mf}(\widehat \bw_i)$, $i=1,2$.
Figure~\ref{fig:ex1_roothaan} finally shows the functions
\[
w_1 \otimes v_1, \quad v_2 \otimes w_2 \quad \text{and}\quad v_1 \otimes v_2
\]
which serve as approximations to the eigenfunctions from Figure~\ref{fig:ex1_full}.  The qualitative structure of the eigenfunctions of the full spatial operator is quite well captured.  

\subsection{Example: a model of n-butane}\label{subsec:butane}

\begin{figure}[htb]
	\centering
		\includegraphics[width=0.35\textwidth]{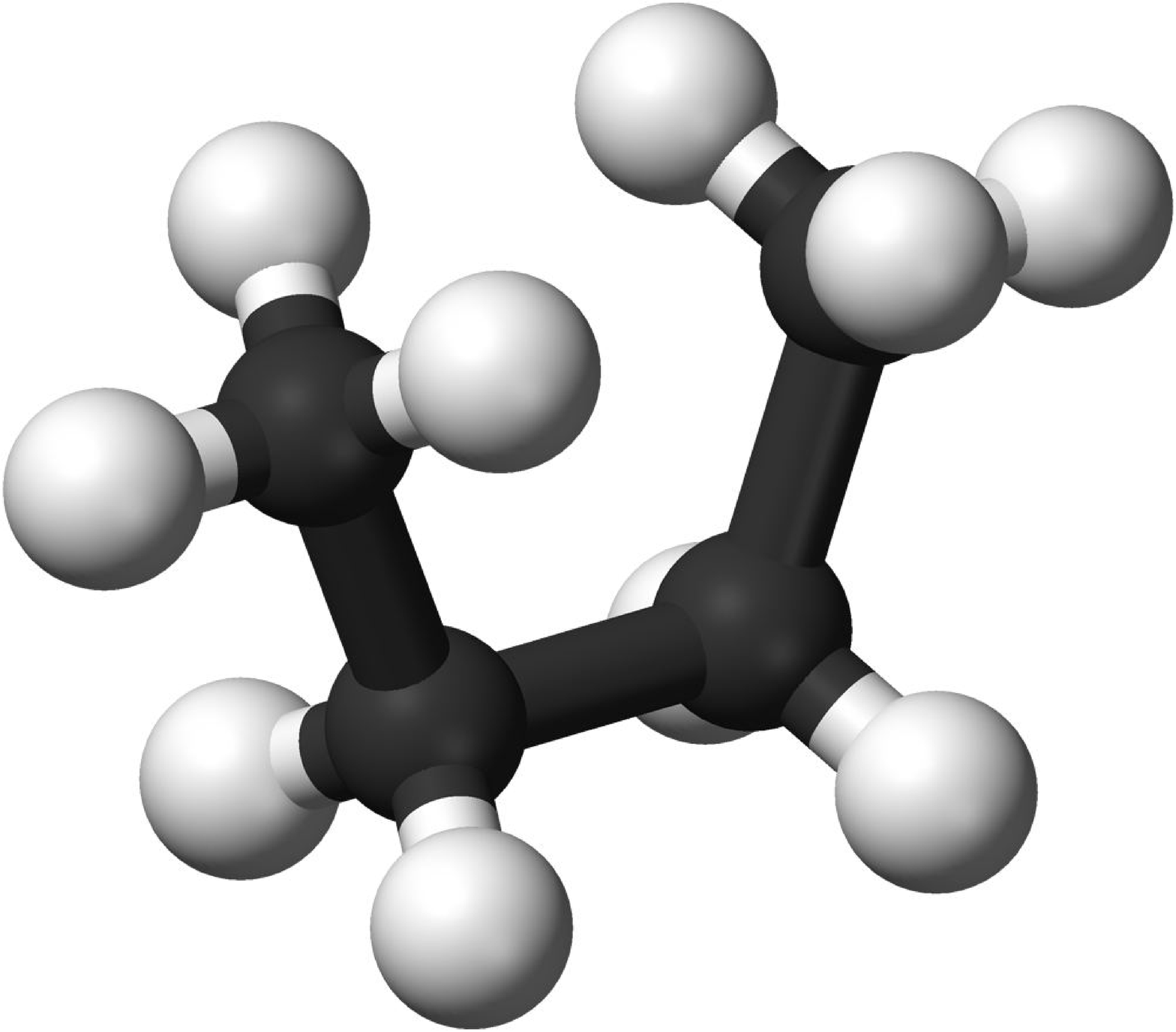}
\qquad
		\includegraphics[width=0.35\textwidth]{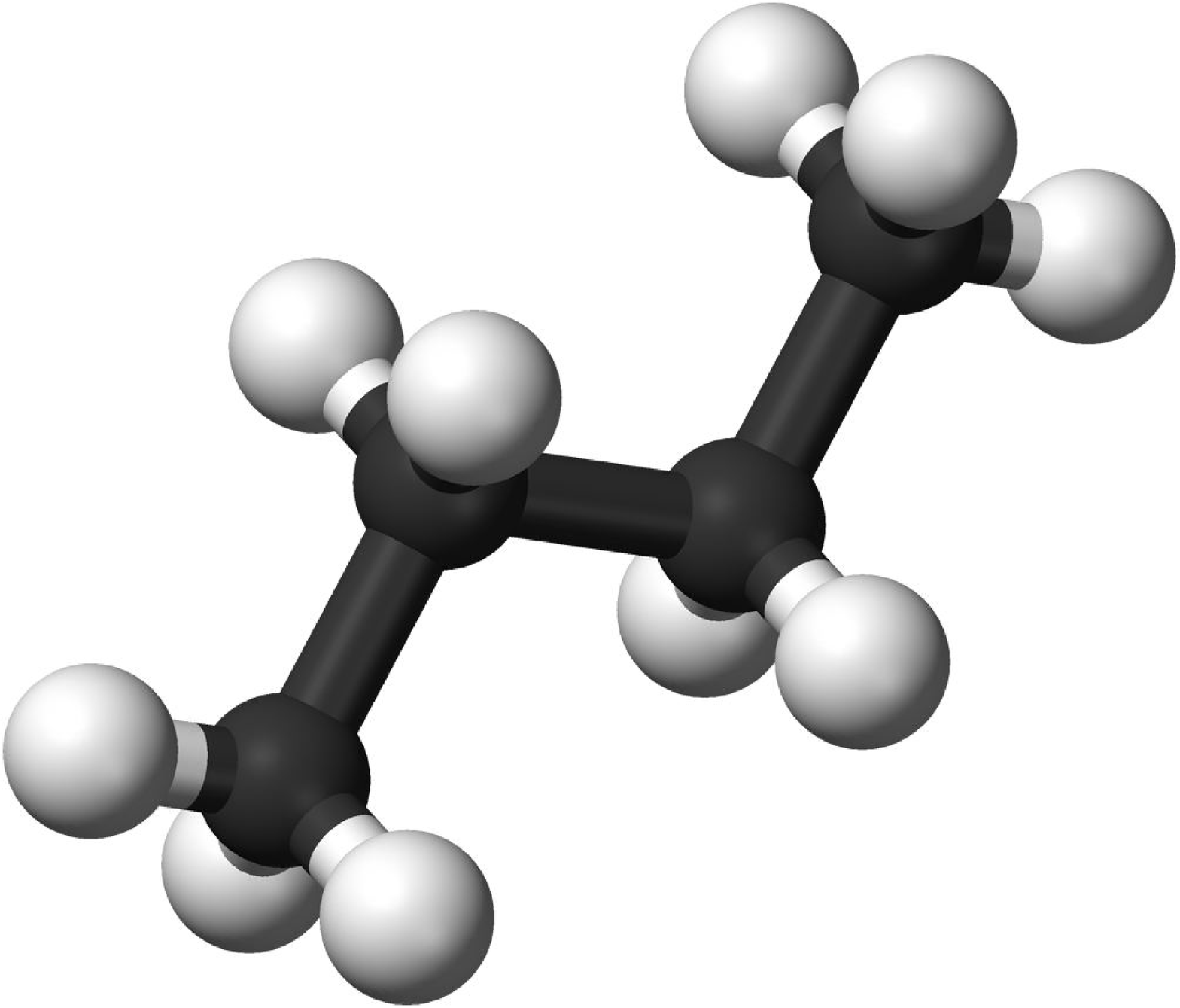}		\caption{Cis- and trans- configuration of $n$-butane.}
	\label{fig:n-butane}
\end{figure}

\begin{figure}[htb]
	\centering
		\includegraphics[width=0.4\textwidth]{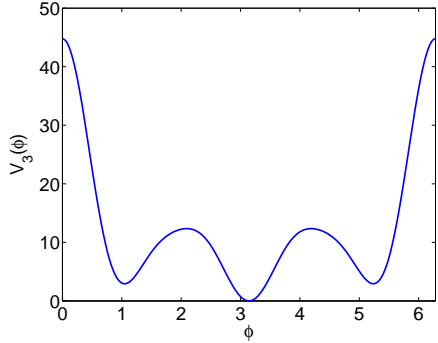}
		\caption{Potential of the torsion angle.}
	\label{fig:torspot}
\end{figure}

As a more realistic test case we analyse the $n$-butane molecule CH$_3$-CH$_2$-CH$_2$-CH$_3$, cf.~Figure~\ref{fig:n-butane}.
More precisely, we consider a united atom model \cite{BrBrOl83a} of this molecule, viewing each CH$_3$ resp.\ CH$_2$ group as a single particle.  Consequently, the configuration of the model is described by six degrees of freedom: three bond lengths, two bond and one torsion angle.  In order to be able to compare the results of our mean field approach to a transfer operator based conformational analysis on the full configuration space, we further simplify the model be fixing the bond lengths at their equilibrium $r_0 = 0.153$ nm. For the bond angles we use the potential
\begin{equation}
		V_{2}(\theta) = -k_{\theta}\left(\cos\left(\theta-\theta_0\right)-1\right)
	\label{eq:anglepot}
\end{equation}
with $k_{\theta}=65 \frac{\text{kJ}}{\text{mol}}$ and $\theta_0=109.47\ ^{\circ}$ and for the torsion angle we employ 
\begin{eqnarray}
		V_{3}(\phi) &  = & K_{\phi}\left(1.116-1.462\cos	\phi-1.578\cos^2\phi+0.368\cos^3\phi\right. \\ \nonumber
		& & \left.+3.156\cos^4\phi+3.788\cos^5\phi\right), \nonumber
	\label{eq:torsionpot}
\end{eqnarray}
with $K_{\phi}=8.314 \frac{\text{kJ}}{\text{mol}}$, cf.~Figure~\ref{fig:torspot}, see also \cite{GrKnZu07a}.
We fix $m_p = 1.672\cdot 10^{-24}$g as the mass of a proton and correspondingly $m_1 = 14\; m_p$ and $m_2 = 15\; m_p$ as the masses of a CH$_2$ and CH$_3$ group, respectively. With $q = (\theta_1,\theta_2,\phi)^{\top}\in [0,\pi]\times[0,\pi]\times[0,2\pi]$ denoting the configuration of our model, the motion of our system is determined by the Hamiltonian
\begin{equation}\label{eq:hamiltonian-butane}
	H(q,p) = \frac{1}{2}p^{\top}M(q)^{-1}p+V(q),
\end{equation}
with $V(q) = V_2(q_1)+V_2(q_2)+V_3(q_3)$ and the mass matrix $M(q)$.  The latter is computed by means of a coordinate transformation $q \mapsto \tilde q(q)$ to cartesian coordinates $\tilde q\in\R^{12}$ for the individual particles, assuming that there is no external influence on the molecule and its linear and angular momentum are zero:  We have
	\[ \dot{\tilde{q}} = D\tilde{q}(q)\dot{q} \]
	and consequently
	\[
	M(q) = D\tilde{q}(q)^{\top}MD\tilde{q}(q),
	\]
	where $M$ denotes the (constant, diagonal) mass matrix of the Hamiltonian in cartesian coordinates.

\subsubsection*{Results for the full operator}

\begin{figure}[tbp]
	\centering
		\includegraphics[width=0.48\textwidth]{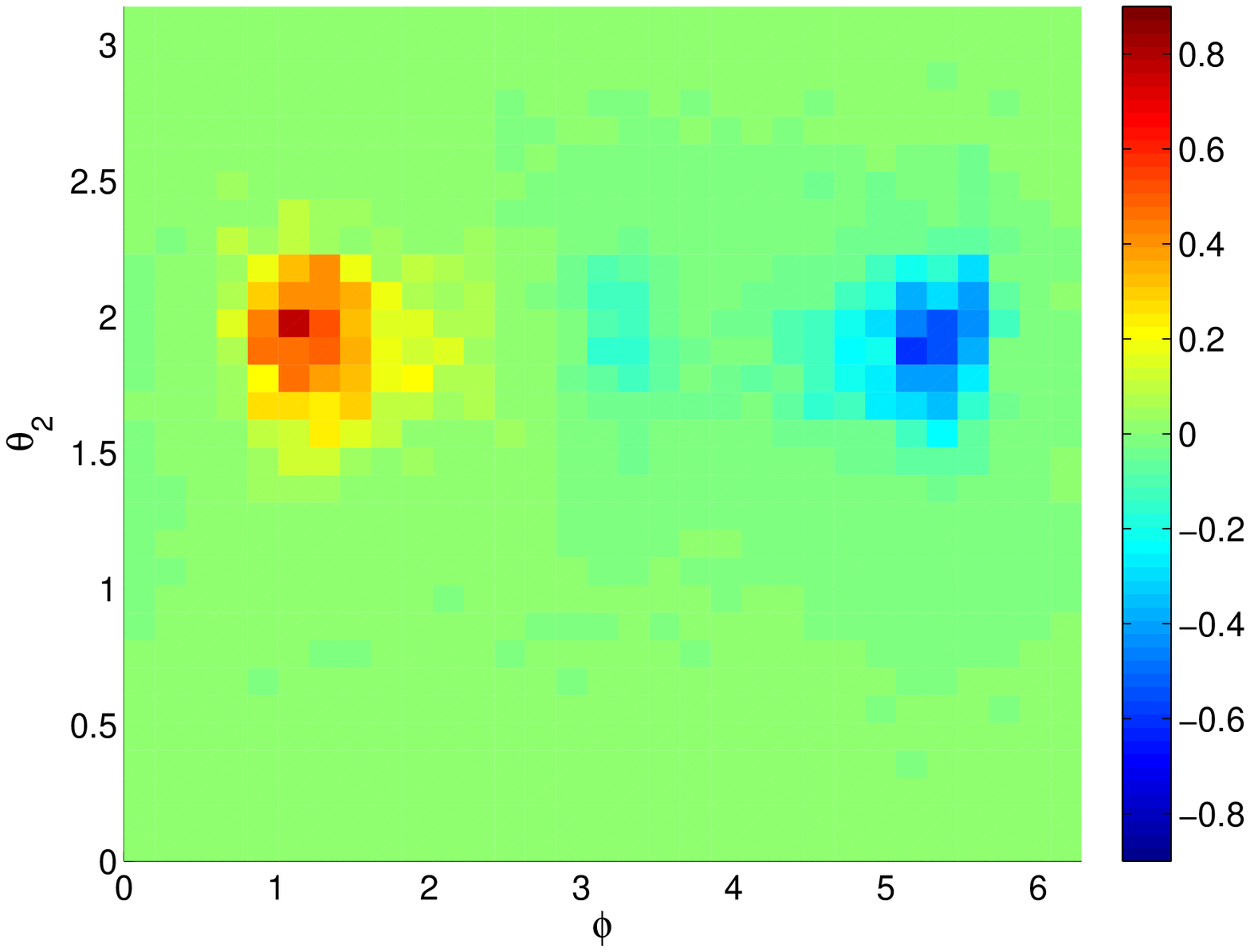}
		\includegraphics[width=0.48\textwidth]{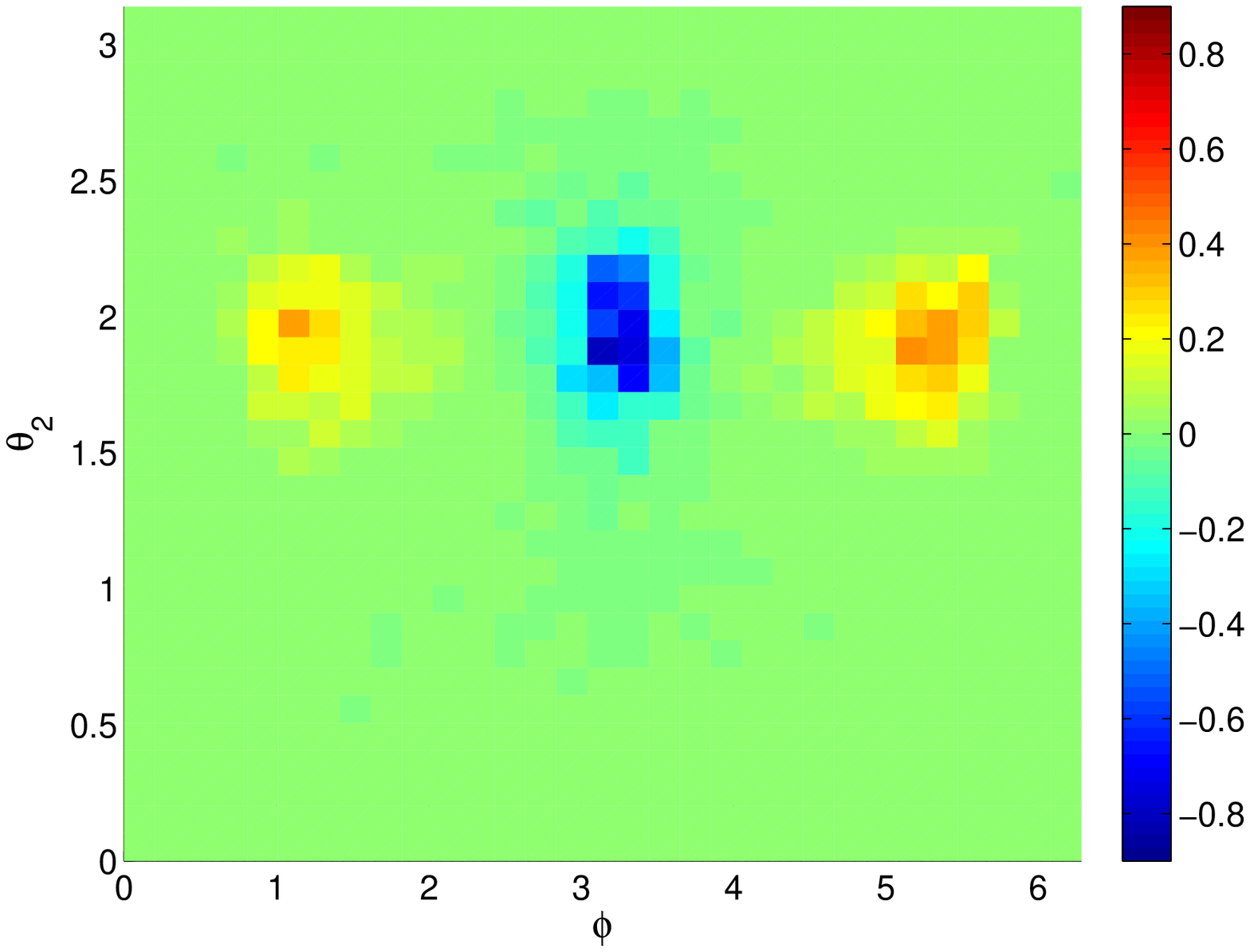}
	\caption{Eigenvectors of the full spatial transfer operator at $\lambda_2=0.985$ (left) and $\lambda_3 = 0.982$ (right). Shown is a slice at $q_1=\theta_1=\pi/2$.}
	\label{fig:butane-evs-full}
\end{figure}

Figure~\ref{fig:butane-evs-full} shows the eigenvectors of the full spatial transfer operator for system (\ref{eq:hamiltonian-butane}) at the two largest eigenvalues $\neq 1$ (computed on a $32\times 32\times 32$ grid using $32$ sample points,  $T=0.5\cdot10^{-13}$ s integration time,  realized by $10$ steps of the explicit Euler method with step size $T/10$).  

\subsubsection*{Results for the mean field approach}

\begin{figure}[btp]
	\centering
		\includegraphics[width=0.48\textwidth]{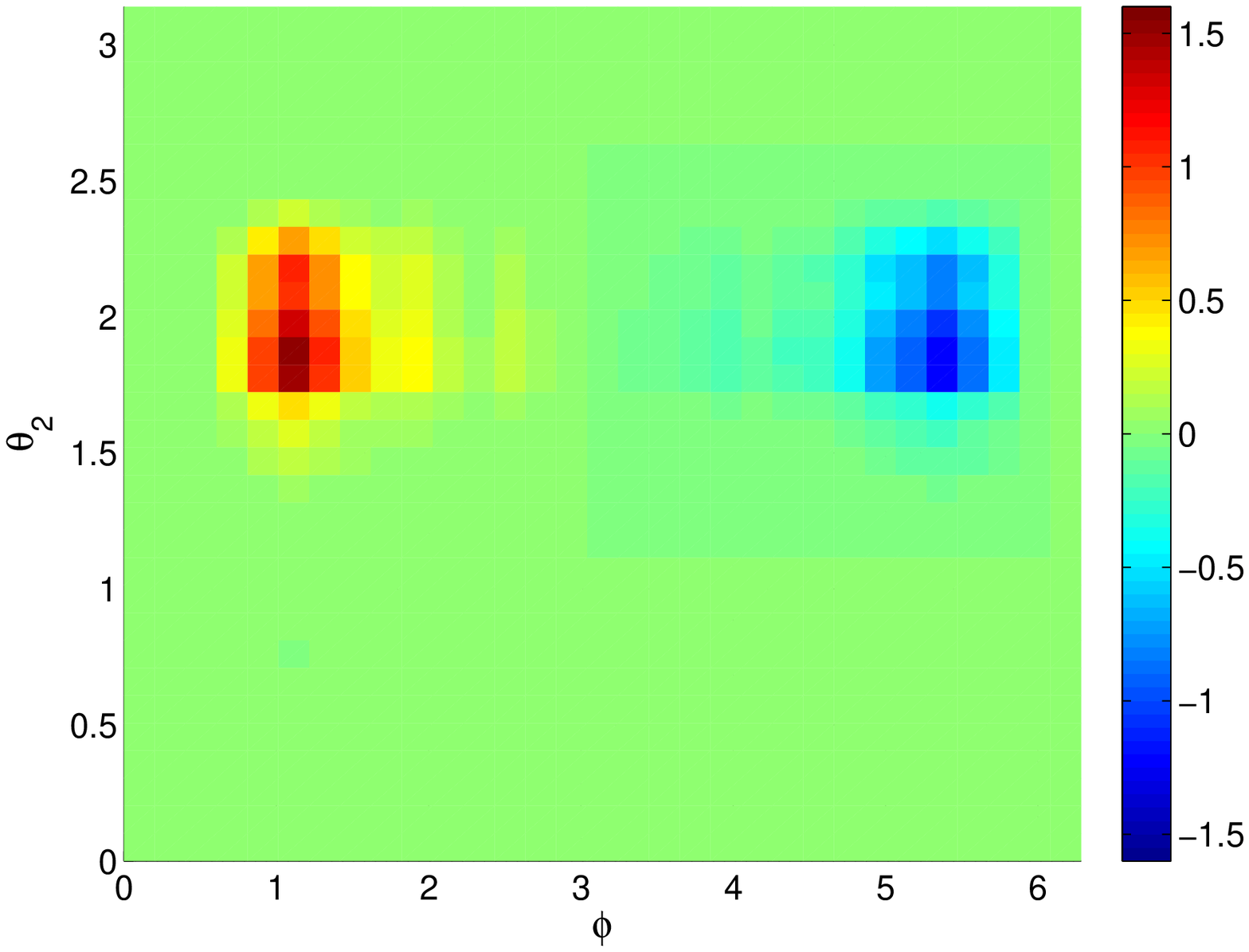}
		\includegraphics[width=0.48\textwidth]{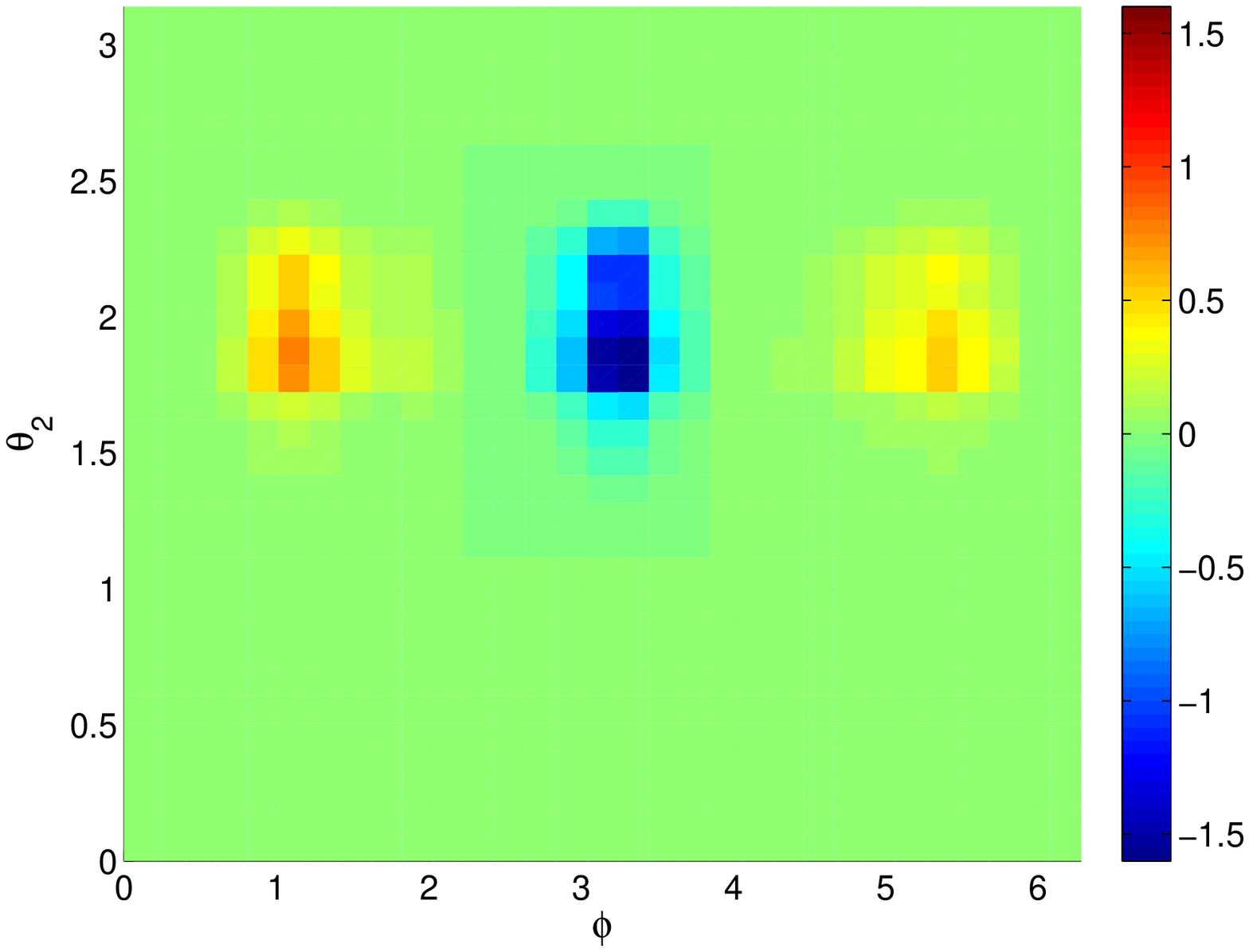}
	\caption{Mean field approximations to the two eigenvectors shown in Fig.~\ref{fig:butane-evs-full}.}
	\label{fig:butane-evs-mf}
\end{figure}

We decompose the model into three subsystems, i.e.\ each configuration variable is treated separately.  The Roothaan iteration is initialized with $w^0_i(q_i):=C_ie^{-\beta V_i(q_i)}$, $i=1,2,3$, where $\beta$ is the inverse temperature corresponding to 300 K and $C_i$ is a corresponding normalizing factor.  Figure~\ref{fig:butane-evs-mf} shows the mean field approximations to the two eigenvectors in Fig.~\ref{fig:butane-evs-full}, namely the products
\[
w_{\theta_2}\otimes v_{\phi,2} \quad\text{and}\quad w_{\theta_2}\otimes v_{\phi,3},
\]
where $w_{\theta_2}$ is the ${\theta_2}$-factor of the invariant density of the mean field system and $v_{\phi,2},v_{\phi,3}$ are the eigenvectors at the second and third largest eigenvalue of the $\phi$ subsystem, respectively (after ten iterations of the Roothan type iteration).
Clearly, the qualitative structure of the eigenvectors of the full operator (cf.~Fig.~\ref{fig:butane-evs-full}) is well captured by the mean field approximation.

\section{Conclusion and Outlook}

We introduced a new approach to the transfer operator based conformational analysis of molecules in this paper.  The central idea is a mean field approach based on a statistical independence ansatz for the eigenfunctions of the operator.  The principal motivation for such an ansatz lies in its linear, as opposed
to exponential, scaling of the number of computational degrees of freedom with the number of atoms, provided the
subsystem size stays fixed.
We established basic theoretical properties including mass and energy conservation and asymptotic correctness
in the limit of weak subsystem coupling. In numerical tests on small systems, the mean field model
is seen to provide a remarkably accurate representation of the true eigenfunctions. Applications to larger systems are currently in progress and will be discussed elsewhere.

\section{Acknowledgements}

We would like to thank Weinan E, Max Gunzburger, Mitch Luskin, and Rich Lehoucq for organizing a very stimulating
Symposium on mathematical issues at the Multiscale materials modeling conference in Tallahassee, Florida, 27-31 October
2008, and for inviting us to present the work described here.
\bibliography{/Users/junge/publications/bibdesk}
\bibliographystyle{abbrv}

\end{document}